\documentstyle[amssymb,10pt]{amsart}
\newtheorem{theorem}{Theorem}[section]
\newtheorem{lemma}[theorem]{Lemma}
\newtheorem{proposition}[theorem]{Proposition}
 
\theoremstyle{definition}  
\newtheorem{definition}[theorem]{Definition}

\newtheorem{conjecture}[theorem]{Conjecture}

\newcommand{\id}{\text{id}}

\newcommand{\nc}{\newcommand}
\nc{\Symm}{{\on{Sym}}}

\newcommand{\on}{\operatorname}   
\newcommand{\eps}{\varepsilon}
 \nc{\cE}{{\cal E}}

\nc{\D}{{\mathfrak d}}
\newcommand{\G}{{\mathfrak{g}}}\nc{\HH}{{\mathfrak h}}
\newcommand{\g}{{\mathfrak{g}}}
\newcommand{\h}{{\mathfrak{h}}}

\nc{\wh}{\widehat}\nc{\wt}{\widetilde}

\newcommand{\la}{{\lambda}}

\newcommand{\ben}{\begin{enumerate}}
\newcommand{\een}{\end{enumerate}}
\newcommand{\ad}{{\text{Ad}}}

\newcommand{\CC}{{\mathbb{C}}}
\newcommand{\QQ}{{\mathbb{Q}}}

\newcommand{\al}{\alpha}
\newcommand{\ZZ}{{\mathbb{Z}}}

\begin{document}

\title{Quantization of Alekseev-Meinrenken dynamical $r$-matrices}

\begin{abstract} The Alekseev-Meinrenken solution $\rho_{\on{AM}}$ to the 
classical dynamical Yang-Baxter equation, associated to a pair 
$(\g,t)$ of a finite dimensional 
Lie algebra $\g$ and $t\in S^2(\g)^\g$, can be 
generalized to an $r$-matrix $\rho$ for a pair $(\g,Z)$, 
where $Z\in \wedge^3(\g)^\g$. 
We construct a dynamical twist $J$ with nonabelian base 
(in the sense of P. Xu), whose classical limit is $\rho$. 
This twist provides a quantization of  
the quasi-Poisson manifold and Poisson groupoid associated 
to $\rho$, and in the case of $(\g,t)$, gives rise to a dynamical 
quantum $R$-matrix, whose classical limit is $r := 
\rho_{\on{AM}} + {t\over 2}$. $J$ is constructed by 
renormalizing associators. We study its analytic
properties in the case of the Knizhnik-Zamolodchikov associator for $\g$, 
introduced by Drinfeld.   
\end{abstract}

\author{Benjamin Enriquez}
\address{IRMA (CNRS), rue Ren\'e Descartes, F-67084 Strasbourg, France}
\email{enriquez@@math.u-strasbg.fr}

\author{Pavel Etingof}
\address{Department of Mathematics, Massachusetts Institute of Technology,
Cambridge, MA 02139, USA}
\email{etingof@@math.mit.edu}

\dedicatory{To the memory of Fridrich Izrailevich Karpelevich}

\maketitle

\section{Introduction}

In \cite{AM}, Alekseev and Meinrenken introduced a 
classical dynamical $r$-matrix $r_{\on{AM}}$, associated 
to a pair $(\g,t)$, where $\g$ is a finite dimensional 
Lie algebra and $t\in S^2(\g)^\g$. 
It turns out that the corresponding modified $r$-matrix 
$\rho_{\on{AM}} := r_{\on{AM}} - {t\over 2}$ depends only on 
$Z = {1\over 4} [t^{1,2},t^{2,3}] \in \wedge^3(\g)^\g$, and that 
the construction of $\rho_{\on{AM}}$ can be generalized to 
a dynamical $r$-matrix $\rho_Z$ for a   
pair $(\g,Z)$, with $Z\in \wedge^3(\g)^\g$
(Section \ref{sect:AM} and \cite{A}). 

In \cite{Xu}, P. Xu defined the notion of quantization of classical
dynamical $r$-matrices with nonabelian base.
In this paper, we provide the first nontrivial example 
of such quantization. 
Namely, following a suggestion of \cite{Xu} (Section 5),
we quantize the dynamical $r$-matrices $\rho_Z$. 

More precisely, we construct a dynamical twist with
nonabelian base in the sense of P. Xu, $J(\lambda)$, such that 
$J(\lambda)=1+\hbar j(\lambda)+O(\hbar^2)$, and
$j-j^{2,1}=\rho_Z$. In the case of $(\g,t)$, 
the quantum dynamical $R$-matrix quantizing 
$r_{\on{AM}}$ is then given by the formula 
$R(\lambda)=J^{2,1}(\lambda)^{-1} * e^{\hbar t/2} * J(\lambda)$, 
where $*$ is a suitable PBW star-product. 
Moreover, in the general case of $(\g,Z)$, 
the dynamical twist $J$ defines a quantization 
of the quasi-Poisson manifold and Poisson groupoid associated to
$\rho_{\on{AM}}$. 

The dynamical twist $J$ is constructed by renormalizing 
(the inverses of) admissible associators attached to $Z$, 
which were constructed in \cite{EnH,Dr2}. 
Namely, let $\Phi \in U(\g)^{\otimes 3}[[\hbar]]$ be such an associator. 
Let us identify the third component $U(\g)$
of this tensor cube with $\Bbb C[\g^*]$ via the
symmetrization isomorphism $S^\cdot(\g)\to U(\g)$, 
and using this identification view 
$\Phi^{-1}$ as a function from $\g^*$ to $U(\g)^{\otimes
2}[[\hbar]]$. We show that there exists a unique series 
$J(\lambda)=1+\hbar j(\lambda)+O(\hbar^2)$, whose 
coefficients are formal functions on $\g^*$ with values in 
$U(\g)^{\otimes 2}$,
regular at $0$, such that $J(\hbar\lambda)=\Phi^{-1}$. 
This element $J(\lambda)$ is the desired dynamical twist. 

We then study the analytic properties of $J(\lambda)$ in the case
of $(\g,t)$ when $\Phi$ given by Drinfeld's Knizhnik-Zamolodchikov (KZ)
associator \cite{Dr2}. We show that $J$ is meromorphic in 
$\lambda$ and regular on an explicitly specified neighborhood of zero. 

We note that the results of this paper generalize in a straightforward 
manner to the case 
when the Lie algebra $\g$ is not necessarily finite dimensional. 
The reason for this is that the subspace $\g_Z$ of $\g$ spanned 
by the first component of $Z$ is a finite dimnensional ideal
(in particular, Lie subalgebra), and in all the arguments $\g$ can be 
replaced with $\g_Z$. 

\medskip\noindent
{\bf Acknowledgements.}
The research of P.E. was partially supported by
the NSF grant DMS-9988796, and was done in part for the Clay Mathematics 
Institute. The authors thank A. Alekseev for a useful
discussion and L. Feh\'er for valuable references. 
P.E. is grateful to V. Toledano Laredo for 
a useful discussion of the dynamical twist equation
with an associator. 

\section{Classical dynamical $r$-matrices
and quasi-Poisson manifolds}

Throughout this paper $\g$ will denote a Lie algebra over $\CC$
with an element $Z \in \wedge^3(\g)^\g$. Unless otherise specified, $\g$ is assumed finite dimensional. We will sometimes 
additionally assume that
$\g$ is equipped with an element $t \in S^2(\g)^\g$, such that 
$Z :=\frac{1}{4}[t^{1,2},t^{2,3}]\in \wedge^3(\g)^\g$.

\subsection{Classical dynamical $r$-matrices}

Let $\h$ be a Lie subalgebra of $\g$. 

\begin{definition} \cite{Fe,EtV1} 
{\it 
A {\rm modified classical dynamical $r$-matrix for} $(\g,Z)$
is a meromorphic $\h$-equivariant map $\rho : \h^*\to 
\wedge^2(\g)$, which satisfies the {\rm modified classical dynamical 
Yang-Baxter equation} (CDYBE) \cite{BDF,Fe}, 
$$
-{\rm Alt}(d\rho)+ \on{CYB}(\rho) = Z,
$$
where $\on{CYB}(\rho) :=  [\rho^{1,2},\rho^{1,3}]+[\rho^{1,2},\rho^{2,3}]+
[\rho^{1,3},\rho^{2,3}]$, and 
$$
{\rm Alt}(d\rho):=\sum h_i^1\frac{\partial \rho^{2,3}}{\partial \lambda_i}-
\sum h_i^2\frac{\partial \rho^{1,3}}{\partial \lambda_i}+
\sum h_i^3\frac{\partial \rho^{1,2}}{\partial \lambda_i}.
$$
Here $(h_i),(\lambda_i)$ are dual bases of $\h$ and $\h^*$.} 
\end{definition}

\begin{definition} 
{\it Assume that $\g$ is equipped with $t\in S^2(\g)^\g$. 
A (quasitriangular) 
{\rm classical dynamical $r$-matrix for} $(\g,t)$  
is a meromorphic $\h$-equivariant function $r: \h^*\to \g\otimes \g$
such that $r(\lambda)+r^{2,1}(\lambda)=t$, which satisfies
the CDYBE 
$$
-{\rm Alt}(dr)+ \on{CYB}(r) =0. 
$$
}
\end{definition}

When $\g$ is equipped with $t \in S^2(\g)^\g$ and $Z := 
{1\over 4} [t^{1,2},t^{2,3}]$, the assignment $r(\lambda)
\mapsto \rho(\lambda) := r(\lambda) - {t\over 2}$ is a 
bijective correspondence between solutions of 
(a) the modified CDYBE for $(\g,Z)$ and (b) the 
CDYBE for $(\g,t)$. 
$\rho(\lambda)$ will be called the modified classical dynamical $r$-matrix 
corresponding to $r(\lambda)$. 

\subsection{Gauge transformations of dynamical $r$-matrices}

Recall from \cite{EtS}, Section 2, that given a dynamical $r$-matrix $\rho$
and a meromorphic function $g: \h^*\to G$ (where $G$ is a 
Lie group with Lie algebra 
$\g$), we can construct another dynamical $r$-matrix
$$
\rho^g:=(g\otimes g)(\rho-\bar\eta_g+\bar\eta_g^{21}-\tau_g)
(g^{-1}\otimes g^{-1}),
$$
where $\eta_g=g^{-1}dg: \in \Omega^1(\h^*)\otimes \g$,
$\bar\eta_g$ is the same element considered as a function
$\h^*\to \h\otimes \g$, and
$$\tau_g(\lambda)=(\lambda\otimes 1\otimes 1)
([\bar\eta_g^{12}(\lambda),\bar\eta_g^{13}(\lambda)]).
$$
This new $r$-matrix is called the gauge transformation of
$\rho$ by $g$. 

\subsection{Formal dynamical $r$-matrices}

It is useful to consider {\it formal} classical dynamical $r$-matrices, i.e., 
those given by formal power series rather than meromorphic
functions. 

The space of formal functions on $\h^*$ is defined as 
$\widehat S^\cdot(\h)$, 
the formal completion of the symmetric algebra 
$S^\cdot(\h)$. 
The formal analogue of the de Rham differential is 
$d : \wh S^\cdot(\h) \to \wh S^\cdot(\h)
\otimes \h$. 

\begin{definition} (\cite{EtS})
{\it A 
formal modified classical dynamical $r$-matrix over $(\g,Z)$
is an $\h$-invariant element of $\widehat S^\cdot(\h) \otimes 
\wedge^2(\g)$, satisfying the modified CDYBE for 
$(\g,Z)$. }
\end{definition}

In the same way, one defines formal quasitriangular dynamical 
$r$-matrices over $(\g,t)$. Formal gauge transformations
are induced by elements $\on{log}(g) \in \wh S^\cdot(\h)_{>0} 
\otimes \g$, where $\wh S^\cdot(\h)_{>0}$ is the maximal ideal of 
$\wh S^\cdot(\h)$ (so they correspond to formal functions $g$ such that 
$g(0) = e$). 

Taylor expansion at the origin takes meromorphic $r$-matrices, 
regular at the origin, to formal $r$-matrices. 

{\bf Remark.} We note that the definitions of a formal $r$-matrix,
formal gauge transformations, etc. 
can be generalized to the case when $\g,\h$ are 
not necessarily finite dimensional. The only changes to be made
are of the following type: e.g., the tensor product  
$\widehat S^\cdot(\h) \otimes 
\wedge^2(\g)$ should be replaced with the completed tensor product 
$\widehat S^\cdot(\h) \widehat\otimes 
\wedge^2(\g)$. 

\subsection{The quasi-Poisson manifold associated to a dynamical $r$-matrix}
\label{sect:23}

Let $\g$ be a finite dimensional Lie algebra 
with a fixed element $Z\in \wedge^3(\g)^\g$. 
Recall \cite{AK,AKM} that a quasi-Poisson $(\g,Z)$-manifold is 
a $\g$-manifold $X$ with an invariant bivector field $\pi$ such that 
the Schouten bracket $[\pi,\pi]$ equals 
$\gamma^{\otimes 3}(Z)$, where $\gamma: \g\to {\rm Vect}(X)$ is the action 
homomorphism\footnote{We will be working with complex 
analytic manifolds; the generalization to 
other categories is straightforward.}. Thus, on a quasi-Poisson manifold 
$X$ we have a skew-symmetric biderivation (quasi-Poisson bracket)
$(f,g)\mapsto\lbrace{f,g\rbrace}$, such that 
$$
\lbrace{\lbrace{f,g\rbrace},h\rbrace}+
\lbrace{\lbrace{g,h\rbrace},f\rbrace}+
\lbrace{\lbrace{h,f\rbrace},g\rbrace}=m(\gamma^{\otimes
3}(Z)(f\otimes g\otimes h)),
$$
where $m$ is the usual multiplication. 

Let $\rho$ be a solution of the modified CDYBE for $(\g,Z)$. 
Let $U\subset \h^*$ be a conjugation-invariant open set 
on which $\rho$ is regular. Using a quasi-Poisson generalization 
of a construction of P. Xu \cite{Xu}, one 
can associate to $\rho$ a quasi-Poisson $(\g,Z)$-structure on the manifold 
$X=U\times G$, where $G$ acts on $X$ by 
$g(u,g')=(u,gg')$. Namely, the quasi-Poisson bracket 
associated to the bivector $\pi$ is given by the formulas 
$$
\lbrace{f_1(u),f_2(u)\rbrace}=
[df_1,df_2](u),
$$
$$
\lbrace{ a,f(g)\rbrace}=(L(a)f)(g),
$$
$$
\lbrace{f_1(g),f_2(g)\rbrace}=(df_1\otimes df_2)(L(\rho(u))),
$$
where $a=a(u)\in \h$ is a linear function, 
and 
for a tensor $x$ in 
$\g^{\otimes k}$, $L(x)$ is the tensor field on $G$ 
obtained by left translations of $x$.

This quasi-Poisson manifold will be denoted by $X_\rho$.

\medskip \noindent
{\bf Remark 1.} Another geometric structure 
which may be attached to a dynamical $r$-matrix is 
the structure of a Poisson groupoid on the groupoid $U\times G\times U$ 
(see \cite{EtV1}).  

\medskip \noindent
{\bf Remark 2.} It is explained in \cite{EtS}, Section 2, that
the Poisson groupoids associated to gauge equivalent dynamical
$r$-matrices are isomorphic. The same argument shows that
the same statement holds for quasi-Poisson manifolds.

\section{Alekseev-Meinrenken $r$-matrices}
\label{sect:AM}

\subsection{The Alekseev-Meinrenken dynamical $r$-matrix for $(\g,t)$}

Let $\g$ be equipped with $t\in S^2(\g)^\g$ and set $\h = \g$. 
Then there exists a unique dynamical $r$-matrix for $(\g,t)$
up to gauge transformations \cite{EtS}. It is  
the Alekseev-Meinrenken dynamical $r$-matrix \cite{AM}. Namely, 
let $t^\vee: \g^*\to \g$ be the linear map associated to $t$.
Then the modified Alekseev-Meinrenken $r$-matrix 
is given by the formula 
$$
\rho_{\on{AM}}(\lambda)=\varphi( 1 \otimes 
{\rm ad}t^\vee(\lambda))(t),
$$
where 
\begin{equation}\label{am}
\varphi(z):=-\frac{1}{z}+\frac{1}{2}{\rm cotanh}\frac{z}{2}.
\end{equation}
$\rho_{\on{AM}}(\lambda)$ takes its values in $\wedge^2(\g)$ 
for the following reason: if $x\in \g$, 
then $\on{ad}(1\otimes x)^k(t)$ lies in $S^2(\g)$ if $k$
is even and in $\wedge^2(\g)$ if $k$ is odd; since $\varphi$
is an odd function, $\rho_{\on{AM}}(\lambda)$ takes its values in 
$\wedge^2(\g)$.

It is clear that $\lambda\mapsto \rho_{\on{AM}}(\lambda)$ is a meromorphic 
modified dynamical $r$-matrix, which is regular (and vanishes) at the origin. 
In particular, it can also be regarded as a formal r-matrix.  

\subsection{The Alekseev-Meinrenken $r$-matrix for a pair $(\G,Z)$}

Let $(\G,Z)$ be a pair of a 
finite dimensional Lie algebra $\G$ and $Z\in \wedge^3(\G)^\G$. 
Set 
$$
Z_n = \Symm_{2,3,\ldots,2n}([Z^{1,2,3},[Z^{3,4,5},\ldots, 
Z^{2n-1,2n,2n+1}]]). 
$$
\begin{lemma}
$Z_n$ is antisymmetric with respect to the exchange of the first and 
last tensor factors. 
\end{lemma}

{\em Proof.} We have 
$Z_n^{2n+1,\ldots,1} = [Z^{2n+1,2n,2n-1},[Z^{2n-1,2n-2,2n-3},\ldots, Z^{3,2,1}]])$. 
Since $Z$ is antisymmetric, this is $(-1)^n 
[Z^{2n-1,2n,2n+1},[Z^{2n-3,2n-2,2n-1},\ldots, Z^{1,2,3}]])$. 
Using $n-1$ times the antisymmetry of the bracket, we find 
$Z_n^{2n+1,\ldots,1} 
= - Z_n$. Then 
$\Symm_{2,3,\ldots,2n}(Z_n^{2n+1,\ldots,1}) 
= \Symm_{2,3,\ldots,2n}(Z_n^{2n+1,2,3,\ldots,2n,1})$ implies the result. 
\hfill \qed\medskip 

$Z_n$ therefore gives
rise to a degree $2n-1$ map $Z^{(n)} : \G^* \to \wedge^2(\G)$, 
defined by $Z^{(n)}(\la) = \langle Z_n, \id\otimes \la\otimes\cdots\otimes 
\la\otimes\id\rangle$. 

Let $(c_n)_{n\geq 0}$ be the numbers such that 
${1\over 2} ( -  {1\over{z}} + \on{cotanh}(z)) 
= \sum_{n\geq 0} c_n z^{2n+1}$. We have 
$c_n = 2^{2n+1} B_{2n+2} / (2n+2)!$, where $B_k$ is the 
$k$th Bernoulli number.

\begin{proposition} 
\label{prop:gen:AM} (also \cite{A})

1) Set $\rho_Z = \sum_{n\geq 0} c_n Z^{(n+1)}$. 
Then $\rho_Z$ is an element of  $(\wh S^\cdot(\G) \otimes \wedge^2(\G))^\G$, 
i.e., it can be viewed as a $\G$-invariant formal map $\G^* \to \wedge^2(\G)$. 
Moreover, it satisfies the modified CDYBE for $(\G,Z)$. 

2) When $\G$ is equipped with an element $t \in S^2(\G)^\G$, and 
$Z = {1\over 4} [t^{1,2},t^{2,3}]$, then $\rho_Z$ coincides with 
$\rho_{\on{AM}}$. 
\end{proposition}

{\em Proof.} Let us first prove 2). Assume that $t\in S^2(\G)^\G$
and set $Z = {1\over 4}[t^{1,2},t^{2,3}]$. Set 
$t_n(\la) = \ad(1\otimes t^\vee(\la))^{2n-1}(t)$ for $\la\in\G^*$. 
As we observed, $t_n$ is a map $\G^*\to \wedge^2(\G)$, polynomial 
of degree $2n-1$. 2) follows from: 

\begin{lemma} 
The maps $t_n$ and $4^n Z^{(n)}$ coincide. 
\end{lemma}

{\em Proof of Lemma.} We have 
$4^n Z_n = [[t^{1,2},t^{2,3}],\ldots,t^{2n,2n+1}]$. 
If we set 
$$
Y_n = [[t^{1,2},t^{2,3}],\ldots,t^{n-1,n}],  
$$
then $4^n Z_n = Y_{2n+1}$. Let us prove by induction 
that for any $n\geq 3$, 
\begin{equation} \label{form:Yn}
Y_n = (-1)^n [[t^{1,n} ,t^{1,n-1}],\ldots,t^{1,2}].  
\end{equation}
This is clear for $n = 3$. Then (\ref{form:Yn}) at step 
$n$ implies that  
\begin{align*}
& Y_{n+1} = [Y_n,t^{n,n+1}] = (-1)^n 
[[[t^{1,n},t^{n,n+1}] t^{1,n-1}],\ldots,t^{1,2}]  
\\ & = 
(-1)^{n+1} [[[t^{1,n+1},t^{1,n}] t^{1,n-1}],\ldots,t^{1,2}] .
\end{align*}
So $4^n Z_n = - [[t^{1,2n+1},t^{1,2n}],\ldots,t^{1,2}]$. Therefore 
\begin{align*}
& 4^n Z^{(n)}(\la) = \langle 4^n Z_n, \id\otimes
\la\otimes\cdots\otimes\la \otimes\id\rangle
\\ & = - [[t,t^\vee(\la)\otimes 1],\ldots,t^\vee(\la)\otimes 1]
= \ad(t^\vee(\la)\otimes 1)^{2n-1}(t) = t_n(\la). 
\end{align*}
\hfill \qed\medskip 

{\em End of proof of Proposition \ref{prop:gen:AM}.}
Let us now prove 1). We are now back to the general case 
$Z\in \wedge^3(\G)^\G$. Clearly, $\rho_Z$ is $\G$-invariant. 
Let us show that it satisfies the modified CYBE. 

Let $\D$ be the Drinfeld 
double of the quasi-Lie bialgebra $(\G,\delta = 0,Z)$. 
$\D$ is equipped with an element $t_\D \in S^2(\D)^\D$. 

Namely, $\D = \G\oplus\G^*$ as a vector space, and its Lie 
bracket is given by $[x,y]_\D = [x,y]$, $[x,\xi]_\D = 
\on{ad}^*(x)(\xi)$, $[\xi,\eta]_\D = \langle Z, \xi\otimes\id\otimes\eta 
\rangle$  for any $x,y\in \G$ and $\xi,\eta\in \G^*$. Here $\on{ad}^*$
is the coadjoint action of $\G$ on $\G^*$. $t_\D$ is dual to the canonical 
pairing of $\G\oplus \G^*$.  We set $Z_\D = {1\over 4} [t_\D^{1,2},t_\D^{2,3}]$. 
Then we have 
\begin{equation} \label{expansion:Z}
Z_\D = Z + \sigma_1 + \sigma_2 + \sigma_3, 
\end{equation}
where $\sigma_1\in \wedge^2(\G^*) \otimes \G$, 
$\sigma_2\in \G^* \otimes \G \otimes \G^*$, and
$\sigma_3\in \G\otimes \wedge^2(\G^*)$. 

Set 
\begin{equation} \label{def:wtZ}
\wt Z_{\D,n} = [Z_\D^{1,2,3},[Z_\D^{3,4,5},\ldots,
Z_\D^{2n-1,2n,2n+1}]] \in \D^{\otimes 2n+1}. 
\end{equation}
Let $\pi : \D \to \G$ be the projection 
on the first summand of $\D = \G\oplus \G^*$. 
We will prove that 
\begin{equation} \label{proj:id}
(\id\otimes \pi^{\otimes 2n-1} \otimes\id)
(\wt Z_{\D,n}) = \wt Z_n + \wt\sigma_{2,n}, 
\end{equation}
where
$$
\wt Z_n = [Z^{1,2,3},[Z^{3,4,5},\ldots,Z^{2n-1,2n,2n+1}]], \; 
\wt\sigma_{2,n} = [\sigma_2^{1,2,3},[\sigma_2^{3,4,5},\ldots,\sigma_2^{2n-1,2n,2n+1}]]. 
$$
Expanding (\ref{def:wtZ}) using (\ref{expansion:Z}), we 
get an expression of $\wt Z_{\D,n}$ with $4^{n}$ terms. 
Applying $\pi$ to the even factors, we retain only
the terms in $Z^{2i-1,2i,2i+1}$ and $\sigma_2^{2i-1,2i,2i+1}$, 
so we have now $2^n$ terms: the 
\begin{equation} \label{expla}
[x_1^{1,2,3},x_2^{3,4,5},\ldots,x_n^{2n-1,2n,2n+1}], 
\end{equation}
where $(x_1,\ldots,x_n) \in \{Z,\sigma_2\}^n$. 
Applying 
$\pi$ to the $(2i+1)$th factor ($2i+1= 3,5,\ldots,2n-1$) 
and recalling that $[\G,\G^*]
\subset \G^*$, we take to zero all terms (\ref{expla})
such that for some $i$, $x_i \neq x_{i+1}$. 
So there remain only two terms:  
the right side of (\ref{proj:id}), 
which is mapped to itself by $\id\otimes \pi^{\otimes 2n-1} \otimes\id$, 
because $[\G^*,\G^*]\subset \G$. 

Recall that $Z^{(n)} = (\Symm_{2,\ldots,2n}(\wt Z_n))^{2n,1,\ldots,2n-1,2n+1}$. 
If we set 
$$
Z^{(n)}_{\D} = (\Symm_{2,\ldots,2n}(\wt Z_{\D,n}))^{2n,1,\ldots,2n-1,2n+1}, 
\; 
\sigma^{(n)}_2 = (\Symm_{2,\ldots,2n}(\wt \sigma_{2,n}))^{2n,1,\ldots,2n-1,2n+1},
$$
and $\pi_{\wh S^\cdot} : \wh S^\cdot(\D) \to \wh S^\cdot(\G)$ the projection associated to 
$\pi$, we get 
\begin{equation} \label{id:2}
(\pi_{\wh S^\cdot} \otimes\id)(Z^{(n)}_\D) = Z^{(n)} + \sigma_2^{(n)}.  
\end{equation}
Let $\rho_{\D}$ be the Alekseev-Meinrenken $r$-matrix attached to 
$(\D,t_\D)$ and set $\tau = (\pi_{\wh S^\cdot} \otimes\id)(\rho_{\D})$. 
Set $\rho' = \sum_{n\geq 0} c_n \sigma_2^{(n+1)}$. Then 
$\tau = \rho_Z + \rho'$, and $\rho' \in \wh S^\cdot(\G) \otimes 
\wedge^2(\G^*)$. 

$\rho_{\D}$ satisfies the CDYBE for $(\G,Z_\D)$: 
$$
-\on{Alt}(d\rho_\D) + \on{CYB}(\rho_{\D})
= Z_\D.  
$$
Apply $\pi_{\wh S^\cdot}\otimes\id$ to both sides. We have 
$(\pi_{\wh S^\cdot} \otimes\id) (d\rho_\D) =
d\tau$ modulo $\wh S^\cdot(\G) \otimes (\G^* \otimes \wedge^2(\D))$. 
Moreover, $(\pi_{\wh S^\cdot} \otimes\id) 
( \on{CYB}(\rho_\D) ) = \on{CYB}(\tau)$. 
Therefore 
\begin{equation} \label{srednye}
-\on{Alt}(d\tau) 
+ \on{CYB}(\tau) = Z_\D \; 
\on{modulo} \; 
\wh S^\cdot(\G) \otimes (\G^* \otimes \wedge^2(\D)).
\end{equation} 

For $i\geq 0$, let $\pi_{\wedge^i} : 
\wedge^i(\D) \to \wedge^i(\G)$ be the 
map induced by the projection $\pi$. 
Apply $\id\otimes \pi_{\wedge^3}$ to (\ref{srednye}).  
Then 
$(\id\otimes \pi_{\wedge^3})(\on{Alt}(d\tau) ) = 
\on{Alt}(d\rho_Z)$. Since $\tau$ has the form 
$\rho_Z + \rho'$, with $\rho' \in \wh S^\cdot(\G) \otimes 
\wedge^2(\G^*)$, we get $(\id\otimes \pi_{\wedge^3})
(\on{CYB}(\tau) ) =\on{CYB}(\rho_Z)$.
Finally, $(\id\otimes \pi_{\wedge^3})(Z_\D) = Z$, and 
$\pi_{\wedge^3}$ takes $\G^* \otimes \wedge^2(\D)$ to zero. 
It follows that $\rho_Z$ satisfies the CDYBE for $(\G,Z)$. 

{\bf Remark.} We note that Proposition \ref{prop:gen:AM} 
remains true if $\g$ is not finite dimnensional.
Indeed, consider the subspace of $\G$ spanned by the 
first components of the tensor $Z$. It is a finite dimensional 
ideal (and hence subalgebra) of $\G$, and $Z\in \wedge^3(\G_Z)^{\G_Z}$. 
The above conclusion then 
holds for the pair $(\G_Z,Z)$. Using the canonical maps 
$\wh S^\cdot(\G_Z) \otimes \wedge^i(\G_Z) \to 
\wh S^\cdot(\G) \otimes \wedge^i(\G)$, we find that 
$\rho_Z$ is a solution of the modified CDYBE for $(\G,Z)$. 
\hfill \qed\medskip 

\begin{proposition} \label{unic:classical}
Let $(\G,Z)$ be a pair of a finite dimensional Lie algebra $\G$ and 
$Z\in \wedge^3(\G)^\G$. Then any formal solution of CDYBE for 
$(\G,Z)$ is gauge-equivalent to $\rho_Z$. 
\end{proposition}

{\em Proof.} In the case $Z=\frac{1}{4}[t^{1,2},t^{2,3}]$,
the result is proved in \cite{EtS} (Theorem 1 and Proposition 3.1).
In the general case, the proof is the same. 
\hfill \qed\medskip

\medskip \noindent
{\bf Remark.} In fact, the formal series  
$\rho_Z$ converges to an analytic function in a neighborhood of zero,
and this function analytically continues to a meromorphic function on 
the whole $\g^*$.
Indeed, (\ref{id:2}) implies that for $\la$
in the formal neighborhood of $0\in \g^*$, 
$\rho_{\D}(\la,0)\in \wedge^2(\g) \oplus \wedge^2(\g^*)$, 
and $\rho_Z(\la)$ is the first component of $\rho_{\D}(\la,0)$.
Since $\rho_{\D}(\la,0)$ converges and extends meromorphically to $\g^*$, 
the same is true for $\rho_Z(\la)$. 


\section{Quantization of classical dynamical $r$-matrices and 
quasi-Poisson manifolds}

In this section, we introduce twist quantization and show its
relation with two problems: (a) the quantization of the 
quasi-Poisson manifold $X_\rho$ (see Section \ref{sect:23}), 
and (b) the quantization of classical dynamical $r$-matrices.  

\subsection{Star-products}

A (local) star-product on a manifold $X$ 
is a product on ${\mathcal O}(X)[[\hbar]]$ 
(where ${\mathcal O}(X)$ is the structure sheaf of $X$),
given by 
$$
f*g=fg+\sum_{m\ge 1}\hbar^mD_m(f,g),
$$
where $D_m$ is a bidifferential operator on $X$ and $\hbar$ is a formal 
variable. 

For example, let $\g$ be a finite dimensional Lie algebra, and $\g^*$
the dual space. Then one can define a star-product on $\g^*$ as
follows. Consider the Lie algebra $\g_\hbar$ over $\Bbb
C[[\hbar]]$, which is $\g[[\hbar]]$ as a vector space, 
with commutator given by the formula $[x,y]_{\g_\hbar}=
\hbar[x,y]_\g$. Consider its universal enveloping algebra
$U(\g_\hbar)$, and identify it with $S^\cdot(\g)[[\hbar]]$ as 
a $\Bbb C[[\hbar]]$-module using the symmetrization isomorphism. 
It is easy to show that this defines a local associative
star-product on $\g^*$ (namely, the product is initially defined on
polynomials, but since it is local, it extends automatically to
the whole structure sheaf). This star-product is called the PBW
star-product, since it is constructed using 
the symmetrization isomorphism, whose bijectivity 
is equivalent to the PBW
theorem. 

\subsection {Quantization of classical dynamical $r$-matrices}
\label{4:2}

According to \cite{Xu}, a
quantization of an $\h$-invariant solution of the 
CDYBE 
$r:\h^*\to \g\otimes \g$ is an $\h$-invariant series 
$$
R \in \on{Mer}(\h^*, U(\g)^{\otimes 2})[[\hbar]],
$$
$R=1+\hbar r+O(\hbar^2)$,  
which satisfies the 
quantum dynamical Yang-Baxter equation (QDYBE)
$$
R^{1,2}(\lambda)*R^{1,3}(\lambda+\hbar h^2)*R^{2,3}(\lambda) =
R^{2,3}(\lambda+\hbar h^1)*R^{1,3}(\lambda)*R^{1,2}(\lambda+\hbar h^3),
$$
where $*$ denotes the PBW star-product of functions on $\h^*$. 
Here we use the dynamical notation $h^{i}$; for example, 
$$
R^{1,2}(\lambda +\hbar h^3):=
\sum_{N\ge 0}
\frac{\hbar^N}{N!}\sum_{i_1,...,i_N=1}^n
(\partial_{\xi^{i_1}}\cdots \partial_{\xi^{i_N}}R)(\lambda)
\otimes (e_{i_1} \cdots e_{i_N})
$$
where $n={\rm dim}(\h)$, and $(e_i)_{i = 1,\ldots,n}$, 
$(\xi^i)_{i = 1,\ldots,n}$ are dual bases of $\h$ and $\h^*$.
This is a
generalization of Felder's equation \cite{Fe} 
from the case when $\h$ is abelian (in which case $*$ 
stands for the usual product). Solutions of QDYBE are called
quantum dynamical $R$-matrices.

\medskip \noindent
{\bf Remark.} If $X$ is a complex manifold and $V$ a complex vector space
(possibly infinite dimensional), then a meromorphic function 
$f: X\to V$ is, by definition, a meromorphic function of $X$ with
values in a finite dimensional subspace of $V$. We denote by 
$\on{Mer}(X,V)$ the space of meromorphic functions $f : X \to V$.

\subsection{Twist quantization} 
\label{4:3}

Let $\G$ be a finite dimensional 
Lie algebra and let $\HH\subset \G$ be a Lie subalgebra. 

\begin{proposition} \label{J-Phi}
Let $\Phi\in U(\G)^{\otimes 3}[[\hbar]]$ and $J\in
\on{Mer}(\h^*, U(\G)^{\otimes 2})[[\hbar]]$ 
be such that: 

(1) $\Phi$ is $\G$-invariant, $\Phi = 1 + O(\hbar^2)$, and 
$\Phi$ satisfies the pentagon equation and the counit axiom (so 
$(U(\G)[[\hbar]],m_0,$ $\Delta_0,\Phi)$ is a quasi-Hopf algebra) 

(2) $J(\la) = 1 + O(\hbar)$ is $\HH$-invariant, and it
satisfies 
\begin{equation} \label{DTE}
J^{12,3}(\la) * J^{1,2}(\la + \hbar h^3) = 
\Phi^{-1} J^{1,23}(\la) * J^{2,3}(\la).  
\end{equation}

Set $Z = \on{Alt}({{\Phi - 1}\over {\hbar^2}})$ mod $\hbar$, 
$j(\la) = ({{J(\la) - 1}\over {\hbar}})$ mod $\hbar$, 
and $\rho(\la) = j(\la) - j(\la)^{2,1}$. Then $Z\in \wedge^3(\G)^\G$, 
and $\rho(\la)$ is a solution of the modified CDYBE for $(\G,Z)$. 
\end{proposition}

{\em Proof.} Straightforward. 
$\square$

When $\Phi$ satisfies the conditions of Proposition \ref{J-Phi}, we
say that it is a quantization of $Z$. When  
$(\Phi,J(\la))$ satisfy the conditions of Proposition \ref{J-Phi}, 
we say that $(\Phi,J(\la))$ is a (twist) quantization 
\footnote{In fact, this definition 
is a generalization of that in \cite{Xu}, 
who considered the triangular case $Z=0$, $\Phi=1$.} 
of the modified CDYBE solution $(Z,\rho(\la))$. 
Equation (\ref{DTE}) is called the 
dynamical twist equation with nonabelian base.

\begin{theorem} \label{Xu}
Let $(\G,t)$ be a pair of a finite dimensional 
Lie algebra $\G$ and $t\in S^2(\G)^\G$. 
Let $\Phi\in U(\G)^{\otimes 3}[[\hbar]]$ be a solution of  
the pentagon and hexagon equations, which satisfies the 
counit axiom (so 
$(U(\G)[[\hbar]],m_0,\Delta_0,e^{\hbar t/2},\Phi)$ is a 
quasitriangular quasi-Hopf algebra), such that
$\Phi$ is a quantization of $Z = {1\over 4}[t^{1,2},t^{2,3}]$. 

Let $r(\la)$ be a solution of the CDYBE for $(\G,t)$. 
Assume that $J(\la)$ is such that $(\Phi,J(\la))$
is a quantization of $(Z,r(\la) - {t\over 2})$. 

Then the function $R(\la) := J^{2,1}(\la)^{-1} * e^{\hbar t/2}
* J(\la)$ is a quantum dynamical $R$-matrix, which quantizes $r(\la)$. 
\end{theorem}

{\em Proof.} For $\h=0$, this was proved by Drinfeld \cite{Dr1}. 
For $\Phi=1$, this is proved in \cite{Xu}. 
The general case is essentially a combination of these two
special cases. \hfill \qed\medskip

\subsection{Formal $R$-matrices and twists}

The definitions and results of Sections \ref{4:2} and \ref{4:3}
extend to the formal case. 
E.g., $J(\la)$ and $R(\la)$ are both elements of 
$(U(\g)^{\otimes 2} \wh\otimes \wh S^\cdot(\h))[[\hbar]]$, 
and $R^{1,2}(\la + \hbar h^3) \in 
(U(\g)^{\otimes 3} \wh\otimes \wh S^\cdot(\h))[[\hbar]]$ 
is the image of $(\id^{\otimes 2} \otimes ((i\otimes \id) 
\circ \Delta))(R)$, where $i : \wh S^\cdot(\h) \to U(\g)[[\hbar]]$ 
is the renormalized symmetrization morphism induced by $\h \ni x\mapsto \hbar x
\in U(\g)[[\hbar]]$ and $\Delta : \wh S^\cdot(\h) \to \wh S^{\cdot}(
\h)^{\wh\otimes 2}$ is induced by $\h\ni x \mapsto x\otimes 1 
+ 1 \otimes x$. 

\medskip \noindent
{\bf Remark.} It is clear that any meromorphic dynamical 
quantum $R$-matrix or twist, which is regular at $0$, 
can be regarded as a formal one by taking the Taylor expansion.  

\begin{conjecture}\label{twistqua} {\it Any 
modified dynamical $r$-matrix admits a (formal) twist quantization. }
\end{conjecture}

\subsection{Quantization of quasi-Poisson manifolds} \label{sect:34}

Let $\g$ be a finite dimensional Lie algebra with a fixed element
$Z\in \wedge^3(\g)^\g$. 
Let $X$ be a quasi-Poisson $(\g,Z)$-manifold, 
and  
$\Phi=1+\frac{\hbar^2}{6}Z+O(\hbar^3)\in (U(\g)^{\otimes
3})^{\g}[[\hbar]]$ 
be an associator for 
$\g$ (so that $(U(\g)[[\hbar]],\Phi)$ is a quasi-Hopf algebra, 
\cite{Dr1}). According to \cite{Dr1}, Proposition 3.10, 
such a $\Phi$ always exists. For example, if
$t\in S^2(\g)^\g$ and $Z=\frac{1}{4}[t^{1,2},t^{2,3}]$, 
$\Phi$ may be taken to be a Lie associator. 

\begin{definition} {\it A quantization of $X$ associated to
$\Phi$ is an invariant 
star-product on $X$, which 
satisfies the equation 
$$
f*g-g*f=\hbar \lbrace{f,g\rbrace}+O(\hbar^2)
$$
(i.e., quantizes the quasi-Poisson bracket on $X$), and 
is associative in the tensor category 
of $(U(\g)[[\hbar]],\Phi)$-modules. This means, 
$$
\mu(\mu\otimes 1)=\mu(1\otimes \mu)\gamma^{\otimes 3}(\Phi), 
$$
on ${\mathcal O}(X)^{\otimes 3}$, 
where $\mu(f\otimes g):=f*g$. }
\end{definition}

\subsection{Quantization of the quasi-Poisson manifold $X_\rho$.}

Let $(\g,Z)$ be a pair of a finite dimensional Lie algebra and
$Z\in \wedge^3(\g)^\g$. Let $\rho$ be a modified classical 
dynamical $r$-matrix for $(\g,Z)$.  
It turns out that to quantize the quasi-Poisson manifold 
$X_\rho$, it suffices to fix a twist quantization of $\rho$. 

Namely, let $\Phi\in (U(\g)^{\otimes 3})^\g[[\hbar]]$ be a Drinfeld associator
for $(\g,Z)$,  $\Phi=1+\frac{\hbar^2}{6}Z+O(\hbar^3)$, 
and let $J$ be a twist quantization of $\rho$ associated to $\Phi$,
which is regular on $U$. 

\begin{theorem} \cite{Xu} \label{quant:poisson}
\footnote{More precisely, this is a generalization of a theorem 
of P. Xu, who considered the special case $Z=0,\Phi=1$; 
the proof, however, is the same as in \cite{Xu}.} 
Consider the star-product on $X_\rho$, which extends the 
PBW star-product on $U$, and satisfies the conditions:

1) Let $f\in U(\h)$, $\Delta(f)=f^{(1)}\otimes f^{(2)}$. Then 
$\psi(g)*f(u)=\psi(g)f(u)$, while
$f(u)*\psi(g)=f^{(2)}(u)(L(f^{(1)})\psi)(g)$, where $L(x)$ is
the differential operator on $G$ obtained by left translations of
$x$; 

2) $\psi(g)*\xi(g)=m(L(J(u))(\psi\otimes \xi))(g)$,
where $m$ is the usual multiplication. 

Then this star-product provides a quantization of the quasi-Poisson 
manifold $X_\rho$. 
\end{theorem}

\medskip \noindent
{\bf Remark.} A twist quantization of $\rho$ also gives rise to a
quantization of the Poisson groupoid corresponding to $\rho$, 
analogously to the construction of \cite{EtV2}. 
To keep the paper short, we will not discuss this construction.

\subsection{Gauge equivalences}

Let $J$ be a formal dynamical twist with a nonabelian base, 
associated to an associator $\Phi$, and $Q
\in (U(\g)\wh \otimes \wh S^\cdot(\h) )^\h[[\hbar]]$, such that 
$(Q$ mod $\hbar)$ has the form $\on{exp}(q_0)$, where 
$q_0 \in (\g \otimes \wh S^\cdot(\h)_{>0})^\h$. 
Set
$$
{}^Q J(\lambda)=Q^{12}(\lambda)*J(\lambda)*Q^2(\lambda)^{-1}*Q^1(\lambda+\hbar
h^2)^{-1}.
$$
It is easy to show that ${}^Q J'$ is also a dynamical twist associated
to $\Phi$. One says that 
${}^Q J$ is {\it gauge equivalent} to $J$ via $Q$. 

Similarly, given a quantum dynamical $R$-matrix $R$, 
define
$$
{}^Q R(\lambda)=Q^2(\lambda+\hbar h^1)*Q^1(\lambda)*R(\lambda)*
Q^2(\lambda)^{-1}*Q^1(\lambda+\hbar
h^2)^{-1}.
$$
It is easy to show
that ${}^Q R$ is also a dynamical $R$-matrix. One says that 
${}^Q R$ is {\it gauge equivalent} to $R$ via $Q$. 

It is clear that the classical limits of gauge equivalent
dynamical twists ($R$-matrices) are gauge equivalent. 
Moreover, the quantized quasi-Poisson manifolds
(Poisson groupoids) associated to gauge equivalent dynamical twists and
$R$-matrices are isomorphic. 

\subsection{Existence of quantization for quasi-Poisson manifolds.}
It is not true that every quasi-Poisson manifold can be quantized. 
Indeed, in the case $Z=0$ the problem of quantization of quasi-Poisson 
manifolds reduces to the problem of equivariant quantization 
of usual Poisson manifolds with a group action preserving 
the Poisson bracket. But it was noted already by Fronsdal in 1978 \cite{Fr}
that such quantization is impossible even in the symplectic case, 
namely for certain coadjoint orbits of semisimple Lie groups. 

However, one may make the following conjecture. 

\begin{conjecture} \label{conj} {\it A quasi-Poisson manifold $X$ 
for $(\g,Z)$ can be 
quantized if the $\g$-action on $X$ is free, i.e., $X$ 
is a principal $G$-bundle over $Y=X/G$, where 
$G$ is the simply connected Lie group corresponding to $\g$.}
\end{conjecture}
 
In the case $G=\{e\}$, this is the Kontsevich formality theorem. 
Moreover, even when $G\ne \{e\}$ but $Z=0$, this conjecture 
apparently follows from formality theory (D.\ Tamarkin and 
V.\ Dolgushev, private communication). 

On the other hand, if $X=G$, then Conjecture \ref{conj} is equivalent to 
the existence of a 
twist quantization of coboundary Lie bialgebras, and is open. 
However, in the quasitriangular case ($X=G$, 
$Z=\frac{1}{4}[t^{1,2},t^{2,3}]$), Conjecture \ref{conj} is equivalent to 
the existence 
of a twist quantization of quasitriangular Lie bialgebras, 
and was proved in \cite{EK}. 

Note also that the quasi-Poisson manifold 
$X_\rho$ attached to a modified dynamical 
$r$-matrix $\rho$ satisfies the freeness assumption, so 
Conjecture \ref{conj} is compatible with Conjecture \ref{twistqua}.

\section{Quantized dynamical twist for the pairs $(\g,Z)$}

In this section, we assume that $(\G,Z)$ is a pair of a 
finite dimensional Lie 
algebra $\G$ and of $Z\in \wedge^3(\G)^\G$. 

\subsection{Admissible associators}

We denote by $U'$ the subspace of $U(\G)[[\hbar]]$ of all 
elements such that $\delta_n(x) \in \hbar^n 
U(\G)^{\otimes n}[[\hbar]]$, 
where $\delta_n = (\id - \eta\circ \varepsilon)^{\otimes n}
\circ \Delta^{(n)}$. $U'$ is the subalgebra generated by the 
$\hbar x$, $x\in \G$. We have $U' /\hbar U' = \wh S^\cdot(\G)$, 
and $U' = U(\G_\hbar)$. 

\begin{definition} (see \cite{EnH})
{\it Let $\Phi\in U(\G)^{\otimes 3}[[\hbar]]$ be a quantization of $Z$. 
$\Phi$ is called {\rm admissible} iff $\hbar \on{log}(\Phi)$
belongs to the $\hbar$-adic completion of 
$(U'_0)^{\otimes 3}$. Here $U'_0$ is the kernel 
of the counit map $U'\to \CC[[\hbar]]$.}   
\end{definition}

In \cite{EnH}, it is proved that there exists an admissible 
quantization $\Phi = \cE(Z)$ of $Z$, defined in terms of $Z$
by universal acyclic formulas (this result is
based on Proposition 3.10 of \cite{Dr1}). 

\subsection{Algebraic versions of dynamical equations} \label{pent}

Let $\h\subset \g$ be a Lie subalgebra. 
Let $J(\lambda)$ be a formal dynamical twist for $\h$ and an associator
$\Phi$ of $\g$, and 
$K(\lambda):=J(\hbar \lambda)$.
This is an $\h$-equivariant function of $\h^*$ with values in 
$U(\g)^{\otimes 2}[[\hbar]]$, polynomial modulo any fixed power 
of $\hbar$, which satisfies the equation 
$$
K^{12,3}(\lambda)\bullet K^{1,2}(\lambda+h^3)=\Phi^{-1}
K^{1,23}(\lambda)\bullet K^{2,3}(\lambda),
$$
where $\bullet$ stands for the multiplication on polynomials 
on $\h^*$ coming from the identification 
$S^\cdot(\h)\to U(\h)$ using the symmetrization map.  

Let us view $K$ as an element of $(U(\g)^{\otimes 2}\otimes U(\h))^\h
[[\hbar]]$. 
(i.e., we view functions of $\lambda$ 
as living the third component $U(\h)$, again using the
symmetrization map). 
Then the twist equation for $K$ becomes an equality 
in $(U(\g)^{\otimes 3}\otimes U(\h))[[\hbar]]$:
\begin{equation}\label{pentform}
K^{12,3,4}K^{1,2,34}=(\Phi^{-1})^{1,2,3}K^{1,23,4}K^{2,3,4}
\end{equation}

Moreover, $K$ has the $\hbar$-adic valuation properties 
\begin{equation} \label{val:K}
K = 1 + O(\hbar), \; 
(\on{id}\otimes \on{id} \otimes \delta_n)(K)
\in \hbar^{n+1} U(\g)^{\otimes 2} \otimes U(\h)[[\hbar]] 
\; \on{if}\; n\geq 1, 
\end{equation}
where $\delta_n : U(\h) \to U(\h)^{\otimes n}$ is $\delta_n = 
(\on{id} - \eta\circ\eps)^{\otimes n} \circ \Delta^{(n)}$, 
$\Delta^{(n)}$ is the $n$fold coproduct of $U(\h)$, 
and $\eps : U(\h) \to \CC$, $\eta : \CC \to U(\h)$ are the counit
and unit maps. 

\begin{proposition}
The map $J(\lambda) \mapsto K$ is a bijection between (a) the set of 
formal solutions of the dynamical twist equations, and (b) the set 
of solutions of (\ref{pentform}), satisfying the valuation conditions 
(\ref{val:K}). 
\end{proposition}

In the same way, if $R(\lambda)$ is a formal solution of QDYBE, 
$S(\lambda) := R(\hbar\lambda)$ and $S\in (U(\g)^{\otimes 2}
\otimes U(\h))^\h [[\hbar]]$ is the image of $S(\lambda)$ by the
symmetrization map, then $S$ satisfies 
\begin{equation} \label{4:QYBE}
S^{1,2,4} S^{1,3,24} S^{2,3,4} = 
S^{2,3,14} S^{1,3,4} S^{1,2,34}
\end{equation}
and the $\hbar$-adic conditions (\ref{val:K}), with $S$ replacing $K$. 

\begin{proposition}
The map $R(\lambda)\mapsto S$ is a bijection between (a) the set of formal 
solutions of QDYBE, and (b) the set of solutions of (\ref{4:QYBE}), 
satisfying conditions (\ref{val:K}). 
\end{proposition}

The map $J(\lambda) \mapsto R(\lambda)$ of Theorem \ref{Xu} 
then corresponds to 
$$
K \mapsto S = (K^{2,1,3})^{-1} e^{\hbar t^{1,2}/2} K^{1,2,3}.
$$ 
The gauge group of all formal functions $Q(\lambda)$ and its actions
also have algebraic counterparts. 

\medskip 

A crucial fact for this paper is that equation (\ref{pentform}) 
strongly resembles the pentagon equation.
We will call it ``the pentagon form of the twist equation''.

\subsection{Quantization of $(Z,\rho_Z)$}

\begin{theorem} \label{thm:Z}
Assume that $\Phi_0\in U(\G)^{\otimes 3}[[\hbar]]$ 
is an admissible quantization of $Z$. Let $J(\la) = 
\Phi_0^{-1}(\hbar^{-1}\la)$, where $\Phi_0^{-1}$ is regarded as
an element of $(U(\G)^{\otimes 2}\otimes S^\cdot(\G))[[\hbar]]$. 

1. The element $J(\la)$ is a formal dynamical twist. More precisely, 
$J(\la) = 1 + \hbar j(\la) + O(\hbar^2) \in (U(\G)^{\otimes 2}
\wh\otimes \wh S^\cdot(\G))[[\hbar]]$ is a series of nonnegative powers 
in $\hbar$ which satisfies the dynamical twist equation. 

2. $J(\la)$ provides a twist quantization of a formal 
CDYBE solution for $(\G,Z)$, gauge-equivalent to $\rho_Z$.   

3. There exists a gauge group element $Q(\la)$, such that 
$J_0(\la) := {}^Q J(\la)$ is a formal dynamical twist, quantizing 
$\rho_Z$. 
\end{theorem}

{\em Proof.} Let us study the $\hbar$-adic valuation properties of 
$J(\la)$. The element $\hbar \on{log}(\Phi_0)$
belongs to the $\hbar$-adic 
completion of $(U'_0)^{\otimes 3}$. Now $U'_0\subset \hbar U(\G)[[\hbar]]$, 
so $\on{log}(\Phi_0) \in \hbar U(\G)^{\otimes 2}[[\hbar]]
\wh\otimes U'_0$. So $\Phi_0 = \on{exp}(\hbar \phi)$, 
where $\phi\in (U(\G)^{\otimes 2} \otimes \wh S^\cdot(\G))[[\hbar]]$. 
this proves that $J(\la)$ has the announced $\hbar$-adic
valuation properties and therefore 1). 

The classical limit of $J(\la)$ is then a solution of the 
modified CDYBE for $(\g,Z)$. According to Proposition 
\ref{unic:classical}, it is gauge-equivalent to $\rho_Z$. 
This proves 2). 

To prove 3), we use the fact that a classical gauge group 
element $g=e^x,\ x \in (\wh S^\cdot(\g)_{>0} \otimes \G)^\G$
may also be regarded as a quantum gauge transformation.  
\hfill \qed\medskip

According to Theorem \ref{quant:poisson}, $J(\la)$ gives rise to 
a quantization of $X_\rho$. 


\subsection{Uniqueness of quantization}

\begin{proposition} \label{5:6}
1. Let $\Phi$ be an associator for $(\g,Z)$ let 
$(\Phi,J_1(\la))$ and $(\Phi,J_2(\la))$ be 
formal twist quantizations of $(Z,\rho_Z)$, then $J_1(\la)$
is gauge-equivalent to $J_2(\la)$.

2. For any associator $\Phi$ for $(\g,Z)$, there
exists a twist quantization $J(\lambda)$ of $\rho_Z$, 
associated to $\Phi$.
\end{proposition}

{\em Proof.}
Proof of 1. As usual, it suffices to show that if $J_1(\lambda)$ and 
$J_2(\lambda)$ coincide 
modulo $\hbar^{n-1}$, then they are gauge equivalent modulo
$\hbar^{n}$. Thus, let us assume that
$$
J_1(\lambda)-J_2(\lambda)=\hbar^n\eta(\lambda)+O(\hbar^{n+1}).
$$ 
The formal function $\eta(\lambda)\in (U(\g)^{\otimes 2}\wh\otimes 
\wh S^\cdot(\g) )^\g$ is a 2-cocycle in the co-Hochschild complex of
$U(\g)$. Thus by adding a $\g$-invariant coboundary (which is an
infinitesimal gauge transformation)
we can assume that $\eta\in (\wedge^2(\g)\wh\otimes
\wh S^\cdot(\g))^\g$ (note that it is not important here whether the
functor of $\g$-invariants is exact or not). If $n=1$, we are
done, as the classical limits of $J_1(\lambda)$ and $J_2(\lambda)$ 
are the same. So let us assume that $n>1$. 

We can think 
of $\eta$ as an invariant formal differential 2-form on $\g^*$. 
Now, it is easy to check that 
the order $\hbar^{n+1}$ of the dynamical twist equation 
tells us that this form is closed (under the de Rham differential).
Thus, by the equivariant Poincar\'e lemma (see, e.g., \cite{EtS}, Lemma
4.1), $\eta=d\xi$, where $\xi$ is an invariant formal 1-form on $\g^*$.
Thus after the gauge transformation $Q(\lambda)=1+\hbar^{n-1}\xi(\lambda)$ 
the twists will be identified modulo $\hbar^{n}$. We are done. 

Proof of 2. According to \cite{Dr1}, Proposition 3.12, 
there exists an
invariant symmetric twist $T\in S^2(U(\g))^\g[[\hbar]]$, which
twists $\Phi$ to the associator ${\cal E}(\widetilde Z)$, where 
$\wt Z\in \wedge^3(\g)^\g[[\hbar]]$ is such that $\wt Z = Z + O(\hbar)$. 
${\cal E}(\wt Z)$ is then admissible. Let $\wt J(\la)$ be the 
twist constructed from ${\cal E}(\wt Z)$ in Theorem \ref{thm:Z}, 3). 
Then $J(\la) := T \wt J(\la)$
is a dynamical twist associated to $\Phi$. 
\hfill \qed \medskip

\section{Quantized dynamical twists for the pairs $(\g,t)$}

In this section, $\g$ is a finite dimensional 
Lie algebra and $t\in S^2(\g)^\g$. 
We set $\h = \g$.

\subsection{}

Let $\g_t$ be the span of the left (or right) tensorands of $t$. 
It is clear that $\g_t$ is an ideal in $\g$. 

\begin{theorem} \label{main}
Assume that $\Phi$ is the image in $U(\g)^{\otimes 3}[[\hbar]]$
of a universal Lie associator. Let 
$J(\lambda):=\Phi^{-1}(\hbar^{-1}\lambda)$,
where $\Phi^{-1}$ is regarded as an element of 
$(U(\g)^{\otimes 2} \otimes \Bbb C[\g^*])[[\hbar]]$. 

1. The element $J$
is a formal dynamical twist. More precisely, 
$J(\lambda)=1+\hbar j(\lambda)+
O(\hbar^2)\in (U(\g)^{\otimes 2}\wh\otimes 
\wh S^\cdot(\g)[[\hbar]]$, 
is a series in nonnegative powers of $\hbar$
which satisfies the dynamical twist equation. 

2. $J(\lambda)$ provides a twist quantization 
of the Alekseev-Meinrenken dynamical $r$-matrix
for $(\g,Z := {1\over 4} [Z^{1,2},Z^{2,3}])$; 
that is, $j-j^{2,1}=\rho_{\on{AM}}$. 

3. The element $R(\lambda):=J^{21}(\lambda)^{-1}*e^{\hbar
t/2}*J(\lambda)$ is a quantum dynamical $R$-matrix, which quantizes 
$r(\lambda) := \rho_{\on{AM}}(\la) + {t\over 2}$. 

4. If $\Phi$ is the Knizhnik-Zamolodchikov associator, 
then $J(\lambda)$ is holomorphic in the open set $U$ of all 
$\lambda\in \g^*$ for which the imaginary part of eigenvalues of 
${\rm ad}t^\vee(\lambda)$ on $\g_t$ belongs to the interval
$(-2\pi,2\pi)$, and extends meromorphically to the whole 
$\g^*$ (in the sense that each 
coefficient of the $\hbar$-adic expansion 
of $J(\lambda)$ is a convergent series
in a neighborhood of zero, and its sum continues to
a meromorphic function on $\g^*$ without poles in $U$). 
\end{theorem}

\subsection{Proof of Theorem \ref{main}}

Proof of 1.  
Let $a=t^{1,2}$, $b=t^{2,3}$. 
Recall \cite{Dr2} that the inverse of a Lie associator is 
a sum of the form $\Phi^{-1}=\sum_w \hbar^{|w|}C_ww(a,b)$, where 
$w(a,b)$ is a word involving letters $a$ and $b$, 
$C_w$ are numbers, $|w|$ is the length of $w$. 
The terms with $w=a^n$ or $b^n$ 
(except $n=0$) do not occur in this expansion
(as any Lie associator equals $1$ modulo the relation $[a,b]=0$). 
When passing from $\Phi^{-1}$ to $J$, for each word $w$ we must use the
symmetrization map to regard the third component of
$w(t^{12},t^{23})$ as an (in general, non-homogeneous) function on
$\g^*$, and divide the part of degree $d$ of this function by
$\hbar^d$. It is clear that $d$ does not exceed the number of
times $b$ occurs in $w(a,b)$, i.e., $d\le |w|-1$ (except for the
empty word). Thus, after division by $\hbar^d$ all powers of
$\hbar$ will be positive except that of the term $1$ 
corresponding to the empty word. Thus, $J=1+\hbar j+O(\hbar^2)$,
as desired. 

That $J$ satisfies the dynamical twist equation follows from
the fact that $\Phi$ satisfies the pentagon relation
$$
\Phi^{1,2,34}\Phi^{12,3,4}=\Phi^{2,3,4}\Phi^{1,23,4}\Phi^{1,2,3},
$$
and subsection \ref{pent}, as equation (\ref{pentform})
has an obvious solution $K=\Phi^{-1}$. 

Another proof of 1. If $\Phi$ is a Lie associator, then 
$\Phi(\hbar t^{1,2},\hbar t^{2,3})$ is an admissible associator.
According to the argument of Theorem \ref{thm:Z}, $\Phi$ 
gives rise to a formal dynamical twist $J(\la)$. 

Proof of 2. The previous argument shows that 
the terms of first order in $\hbar$ are obtained from words $w$ 
in which $a$ occurs only once, by taking the principal symbol of
elements of $U(\g)$ that arise in the third component. 
Since $\Phi$ is a Lie associator (the exponential of a Lie
series) the words $w$ of length $l$ with only one occurrence of $a$ combine
into a single term $({\rm ad}b)^{l-1}(a)$. 
Thus, the element
$\tilde\rho(\lambda):=j(\lambda)-j^{2,1}(\lambda)$
has the form
$$
\tilde\rho(\lambda)=\tilde\varphi({\rm ad}t^\vee(\lambda))^{(2)}(t),
$$ 
where $\tilde\varphi\in {\Bbb C}[[z]]$ is a formal series
containing only odd powers of $z$ (a priori, depending on $\Phi$)
and the superscript (2) means the action in the second
component. 
Also, $\tilde\rho$ satisfies the 
modified classical dynamical Yang-Baxter equation.
Thus, taking $\g$ to be any finite dimensional 
simple Lie algebra (e.g., ${{\mathfrak sl}}(2)$) 
and using Theorem 3.14 of \cite{EtV1}
and the classification of dynamical $r$-matrices for $\g$ (\cite{EtV1}), 
we find that the only possibility for $\tilde\varphi$ is
$\tilde\varphi=\varphi$ (for another, direct proof of this fact see 
\cite{FP}, Section 3). 
Thus $\tilde\rho=\rho_{\on{AM}}$, as desired. 

\medskip \noindent
{\bf Remark.} In the course of this proof we have established the
(well known) fact that the coefficients
of $({\rm ad}b)^{2k-1}a$
in $\Phi$ are rational numbers which 
do not depend of $\Phi$ (and express via Bernoulli
numbers). For the KZ associator, this leads to the classical
formula for $\zeta(2k)$.  
\medskip 

Proof of 3. The statement follows from 1,2, and Theorem
\ref{Xu}. 

Proof of 4. 
Recall \cite{Dr2}  
that the KZ associator $\Phi$ 
is the ``monodromy'' from $0$ 
to $1$ of the differential equation 
$$
\frac{dF}{dz}=\frac{\hbar}{2\pi i}\left(\frac{a}{z}+\frac{b}{z-1}\right)F
$$
(where $F(z)\in U(\g)^{\otimes 3}[[\hbar]]$). Therefore, 
$$
\Phi=\lim_{u\to 1-}(1-u)^{-\frac{\hbar b}{2\pi i}}
\left(\sum_{n=0}^\infty \left(\frac{\hbar}{2\pi i}\right)^n
\int_{0<z_1\le...\le
z_n<u}\frac{(\frac{u}{z_n})^{\frac{\hbar a}{2\pi
i}}b(\frac{z_n}{z_{n-1}})^{
\frac{\hbar a}{2\pi
i}}...(\frac{z_2}{z_1})^{\frac{\hbar a}{2\pi i}}bz_1^{\frac{\hbar a}{2\pi i}}}
{(z_n-1)...(z_1-1)}d\mathbf z\right)
$$
Since $\Phi^{-1}=\Phi^{3,2,1}$, it is given by the same formula 
with $a$ and $b$ interchanged. Namely, using the notation ${\rm
ad}(x)(y):=[x,y]$ and changing variables, we have  
$$
\Phi^{-1}=\lim_{u\to 1-}(1-u)^{-\frac{\hbar a}{2\pi i}}
\left(\sum_{n=0}^\infty \left(\frac{\hbar}{2\pi i}\right)^n
\int_{0<z_1\le...\le
z_n<1}\frac{z_n^{-\frac{\hbar {\rm ad}b}{2\pi
i}}(a)...z_1^{-\frac{\hbar {\rm ad} b}{2\pi i}}(a)}
{(uz_n-1)...(uz_1-1)}d\mathbf z\right)
$$
This implies that the element $J$ is given by the formula 
$$
J(\lambda)=
\lim_{u\to 1-}(1-u)^{-\frac{\hbar t}{2\pi i}}
\left(\sum_{n=0}^\infty \left(\frac{\hbar}{2\pi i}\right)^n
\int_{0<z_1\le...\le
z_n<1}\frac{z_n^{-\frac{{\rm ad}t^\vee(\lambda)_*^{(2)}}{2\pi
i}}(t)*...*z_1^{-\frac{{\rm ad} t^\vee(\lambda)_*^{(2)}}{2\pi i}}(t)}
{(uz_n-1)...(uz_1-1)}d\mathbf z\right).
$$
Here the symbol ${\rm ad}x_*$ denotes the operator
of commutation with $x$ in the algebra 
$(U(\g)^{\otimes 2}\wh\otimes \wh S^\cdot(\g)[[\hbar]]$,
where the third factor is equipped with the PBW star-product, 
and the superscript $(2)$, as usual, denotes the action in the
second component.  

We will now need the following simple lemma from complex analysis. 

\begin{lemma} \label{compan} Let $u<1$, and $k_1,...,k_n$ be
nonnegative integers. 
Consider the function 
$$
F_{n,u}^{k_1,...,k_n}(\alpha_1,...,\alpha_n):=
\int_{0<z_1\le...\le
z_n<1}\frac{z_1^{\alpha_1}...z_n^{\alpha_n}\ln^{k_1}z_1...\ln^{k_n}z_n}
{(1-uz_1)...(1-uz_n)}d\bold z,
$$
defined for ${\rm Re}(\alpha_j)>0$. 
Then 

(i) the function $F_{n,u}^{k_1,...,k_n}$ 
extends meromorphically to all values of $\alpha_1,...,\alpha_n$, 
and its poles are on the hyperplanes $\alpha_1+...+\alpha_k+k+m=0$,
$k=1,...,n$, and $m$ is a nonnegative integer, and

(ii) for any $\varepsilon>0$ one has  
$$
F_{n,u}^{k_1,...,k_n}(\alpha_1,...,\alpha_n)=
\sum_{p=0}^n F_n^{p,k_1,...,k_n}(\alpha_1,...,\alpha_n)\ln^p(1-u)
+O((1-u)^{1-\varepsilon})
$$
as $u\to 1-$, where $F_n^{p,k_1,...,k_n}$ are meromorphic with poles 
as in (i).  
\end{lemma} 

{\em Proof.}
Let us prove (i). It suffices to prove the result 
for $k_i=0$; the general case can be recovered by
differentiation in $\alpha_i$.
Expanding in powers of $u$ and computing the integral, we get 
$$
F_{n,u}^{0,...,0}(\alpha_1,...,\alpha_n)=
\sum_{p_1,...,p_n\ge 0}\frac{u^{p_1+...+p_n}}{\prod_{j=1}^n
(\alpha_1+...+\alpha_j+j+p_1+...+p_j)}.
$$
One checks that the radius of convergence of the r.h.s. is $1$ 
when $(\alpha_1,...,\alpha_n)$ does not belong to the hyperplanes, 
which  implies statement (i). 

Let us prove (ii). Set ${\bold k} = (k_1,...,k_n)$. 
$$
G_{n,u}^{{\bold k}}(\al_1,...,\al_n)
= \int_{0 \leq z_1< ...< z_n \leq u}
{ {z_1^{\al_1}  ... z_n^{\al_n} \ln^{k_1}(z_1)  ... \ln^{k_n}(z_n)}
\over {(1-z_1)...(1-z_n)} } d \bold z. 
$$
We have 
$$
F_{n,u}^{{\bold k}}(\al_1,...,\al_n) = u^{-(\al_1 + ... + 
\al_n + n)} \sum_{m,{\bold l}}
\lambda_{n,{\bold k}}^{m,{\bold l}}
(\ln^{m}u) G^{{\bold l}}_{n,u}(\al_1,...,\al_n), 
$$
where the $\lambda_{n,{\bold k}}^{m,{\bold l}}$ are 
constants. We will show that 
$G^{{\bold l}}_{n,u}(\al_1,...,\al_n)$ has an expansion 
$$
G^{{\bold l}}_{n,u}(\al_1,...,\al_n)
= \sum_{p=0}^n G_n^{p,{\bold l}}(u,\al_1,...,\al_n) \ln^p(1-u), 
$$
where $G_n^{p,{\bold l}}(u,\al_1,...,\al_n)$ is a continuous 
function on 
$[0,1] \times \CC^n$, smooth in $u\in (0,1]$ and meromorphic in 
$(\al_1,...,\al_n)$ with poles as in (i). 

Let us treat the case where all Re$(\al_i)$ are $>0$. We use the following
fact. Let $f(z,\underline\al)$ be a continuous function defined on $[0,1]\times 
\{\al\in \CC | \on{Re}(\al)>0\}^n$, smooth in $z\in (0,1]$ and holomorphic in 
$\underline{\al}$ (that is, $f(-,\underline{\alpha})$ is smooth for any 
$\underline{\alpha}$,  and $f(z,-)$ is holomorphic
for any $z$). For any $k\geq 0$, we have 
$$
\int_0^z f(t,\underline\al) {{\ln^k(1-t)}\over{1-t}} dt
= - {{f(1,\underline{\al})}\over{k+1}} \ln^{k+1}(1-z) + 
g(z,\underline{\alpha}), 
$$
where $g(z,\underline{\alpha}) = \int_0^z {{ f(t,\underline{\al}) - f(1,\underline{\al})}
\over{1-t}} \ln^k(1-t)dt$ has the same properties as $f(z,\underline{\alpha})$. 
We apply this result to $\int_{0<z_1<z_2}
{{z_1^{\al_1} \ln^{k_1}z_1}\over{1-z_1}} dz_1$ (with $k = 0$). 
This expression is then of the form $g_1(z_2,\al_1)\ln(1-z_2) 
+ g_2(z_2,\al_1)$. We transport it in the integral and iterate the procedure. 
We then extend the result to $\CC^n$ by induction on $n$, using the 
functional equation relating 
$G_{n,u}^{{\bold k}}(\underline\alpha - \underline\delta_i)
- G_{n,u}^{{\bold k}}(\underline\alpha)$ with the 
$G_{n-1,u}^{{\bold l}}(\underline\alpha + \underline\al')$, 
$\underline\al'\in\ZZ^{n-1}$ (here $\underline\delta_i$ is the 
$i$th basis vector of $\ZZ^n$). 
\hfill \qed \medskip 

Let us now continue the proof of 4. 
Fix $\lambda$ and 
let elements $C_n(u,\lambda)\in U(\g)\otimes U(\g)$ be defined by the formula
$$
\sum_{n\ge 0}C_n(u,\lambda)\hbar^n:=
\left(\sum_{n=0}^\infty \left(\frac{\hbar}{2\pi i}\right)^n
\int_{0<z_1\le...\le
z_n<1}\frac{z_n^{-\frac{{\rm ad}t^\vee(\lambda)_*^{(2)}}{2\pi
i}}(t)*...*z_1^{-\frac{{\rm ad} t^\vee(\lambda)_*^{(2)}}{2\pi i}}(t)}
{(uz_n-1)...(uz_1-1)}d\bold z\right).
$$
Since ${\rm ad}x_*={\rm ad}x+O(\hbar)$, 
the coordinates of $C_n(u,\lambda)$ with respect to any basis of 
$U(\g)\otimes U(\g)$ are linear
combinations with complex coefficients of functions of the form 
$F_{n,u}^{k_1,...,k_n}(\alpha_1,...,\alpha_n)$, 
where $-2\pi i\alpha_j$ are eigenvalues of the operator 
${\rm ad}t^\vee(\lambda)$ (to see this, it suffices to 
bring this operator to the Jordan normal form). This shows that 
$C_n(u,\lambda)$ is well defined and holomorphic in $\lambda$ when $\lambda$ 
belongs to the open set $U$ (in this case the integral 
which defines $F_{n,u}^{k_1,...,k_n}(\alpha_1,...,\alpha_n)$ is 
absolutely convergent). 

Now, we have $J(\lambda)=\sum_{N\ge 0}\hbar^NJ_N(\lambda)$, where 
$J_N(\lambda)$ is given by the formula 
$$
J_N(\lambda)=\lim_{u\to 1-}
\sum_{k=0}^N\frac{(-t/2\pi i)^k}{k!}C_{N-k}(u,\lambda)\ln^k(1-u).
$$
From this and Lemma \ref{compan} (ii), it follows 
that $J_N$ is also holomorphic 
in $U$. The existence of the meromorphic 
extension of $J_N$ to the whole $\g^*$ now follows from Lemma \ref{compan}, 
(i),(ii).
The theorem is proved. 

\medskip \noindent
{\bf Remark.} Let us confirm by an explicit computation that
for the KZ associator, $j-j^{21}=\rho_{\on{AM}}$. 
According to the above formulas, we have 
$$
j(\lambda)=-\frac{1}{2\pi i}\lim_{u\to 1_-}
\left(t\ln(1-u)+\left(\int_0^1\frac{z^{-{\rm
ad}t^\vee(\lambda)/2\pi i}}{1-uz}dz\right)^{(2)}(t)\right).
$$
Computing this explicitly, we get 
$$
j(\lambda)=(2\pi i)^{-1}\left(-\frac{1}{1-\frac{1}{2\pi i}{\rm ad}t^\vee(\lambda)}
+\sum_{p\ge 1}\left(\frac{1}{p}-
\frac{1}{p+1-\frac{1}{2\pi i}{\rm ad}t^\vee(\lambda)}\right)\right)^{(2)}(t)
$$
Now, observe that the space $\g_t$ carries a nondegenerate inner product 
$t^{-1}$, with respect to which the operators 
${\rm ad} x$, $x\in \g$, are skew-symmetric. 
Hence, 
$$
j(\lambda)-j^{2,1}(\lambda)=
(2\pi i)^{-1}\left(\lim_{N\to \infty}\sum_{p=-N,p\ne 0}^N
\frac{1}{\frac{1}{2\pi i}{\rm ad}t^\vee(\lambda)-p}\right)^{(2)}(t)
$$
Recalling that $\lim_{N\to \infty}\sum_{p=-N}^N \frac{1}{z-p}=\pi
{\rm cotan} \pi z$, we find that
$$ 
j(\lambda)-j^{2,1}(\lambda)=\varphi({\rm ad}t^\vee(\lambda))^{(2)}(t),
$$
where $\varphi$ is defined by formula (\ref{am}).

\subsection{Uniqueness of quantization} 

\begin{proposition} 1. Let $\Phi$ be any associator for $\g$. 
A twist quantization $J(\lambda)$ of $\rho_{\on{AM}}$
associated to $\Phi$ (which exists according to Proposition
\ref{5:6}, 2) can be chosen to be holomorphic
on the open set $U$. 

2. The quantum dynamical $R$-matrix $R(\lambda)$ 
attached to a Lie associator $\Phi$
in Theorem \ref{main} is in fact independent of $\Phi$ up to a
gauge transformation. 
\end{proposition}

{\em Proof.}
Proof of 1. By Drinfeld's results \cite{Dr1,Dr2}, there exists an
invariant symmetric twist $T\in S^2(U(\g))^\g[[\hbar]]$, which
twists $\Phi$ to the Knizhnik-Zamolodchikov associator
$\Phi_{\on{KZ}}$. Now let $J(\lambda)$ be the holomorphic dynamical twist 
constructed in Theorem \ref{main}, and let $J'(\lambda)=TJ(\lambda)$. 
Then $J'(\lambda)$ is the desired dynamical twist associated to $\Phi$.  

Proof of 2. Let $\Phi,\Phi'$ be two Lie associators, 
and $J(\lambda),J'(\lambda)$, $R(\lambda),R'(\lambda)$ 
the corresponding dynamical twists 
and $R$-matrices constructed in Theorem \ref{main}. 
Let $T$ be a Drinfeld twist mapping $\Phi$ to $\Phi'$.
Then $J(\lambda)$ and $TJ'(\lambda)$ satisfy the conditions of 
Proposition \ref{5:6}, 1), 
so they are gauge equivalent. In particular, their associated $R$-matrices
are gauge equivalent. The $R$-matrix of $J(\lambda)$ is $R(\lambda)$, 
and since $T$ is symmetric and invariant, the $R$-matrix of $T J'(\lambda)$
is the same as that of $J'(\lambda)$, i.e., $R'(\lambda)$. So 
$R(\lambda)$ and $R'(\lambda)$ are gauge equivalent.
\hfill \qed 
\end{document}